\scrollmode
\documentstyle[12pt,amssymb]{article}
\def\R{{\Bbb R}}

\def\N{{\Bbb N}}
\def\eas{\begin{eqnarray*}}
\def\eeas{\end{eqnarray*}}
\def\lll{\lefteqn}
\def\nn{\nonumber}
\def\eq#1{\begin{equation}\label{#1}}
\def\eeq{\end{equation}}
\def\ea#1{\begin{eqnarray}\label{#1}}
\def\eea{\end{eqnarray}}
\def\la{\label}
\def\re#1{(\ref{#1})}
\def\dist{\mathop{\mbox{\rm dist}}}

\def\dev{\mathop{\mbox{\rm dev}}}
\def\div{\mathop{\mbox{\rm div}}}
\def\sym{\mathop{\mbox{\rm sym}}}
\def\skew{\mathop{\mbox{\rm skew}}}
\def\div{\mathop{\mbox{\rm div}}}
\def\Curl{\mathop{\mbox{\rm Curl}}}
\def\tr{\mathop{\mbox{\rm tr}}}
\def\id{\mathop{\mbox{\rm id}}}
\def\Sing{\mathop{\mbox{\rm Sing}}}
\def\a{{\alpha}}
\def\b{{\beta}}
\def\g{{\gamma}}
\def\d{\partial}
\def\th{\theta}

\def\l{{\lambda}}
\def\p{{\varphi}}
\def\dx{\,dx}

\def\HH{{\cal H}}
\def\dH{\,d{\cal H}}
\def\del{{\delta}}
\def\eps{{\varepsilon}}

\def\qed{\hfill$\Box$\\\strut}

\def\<{\langle}
\def\>{\rangle}
\newtheorem{lemma}{Lemma}[section]

\newtheorem{theorem}[lemma]{Theorem}
\newtheorem{proposition}[lemma]{Proposition}
\def\bib#1{\bibitem[#1]{#1}}

\def\tf#1#2{{\textstyle\frac{#1}{#2}}}
\def\loc{{\mbox{\tiny\rm loc}}}

\def\tp{\widetilde{\p}}
\def\tR{\widetilde{R}}

\def\opi{\overline{\p}_{(i)}}
\def\oRi{\overline{R}_{(i)}}
\def\drad{\d_{\mbox{\tiny rad}}}
\def\mue{\mu_1}
\def\muc{\mu_c}
\def\muz{\mu_2}

\def\XXint#1#2#3{{\setbox0=\hbox{$#1{#2#3}{\int}$}
                 \vcenter{\hbox{$#2#3$}}\kern-.5\wd0}}
\def\Xint#1{\mathchoice
              {\XXint\displaystyle\textstyle{#1}}
              {\XXint\textstyle\scriptstyle{#1}}
              {\XXint\scriptstyle\scriptscriptstyle{#1}}
              {\XXint\scriptscriptstyle\scriptscriptstyle{#1}}
              \!\int}
\def\mint{\Xint-}

\begin{document}
\sloppy
\title{Regularity issues for Cosserat continua\\ and
$p$-harmonic maps}
\author{Andreas Gastel}
\maketitle

\begin{center}
\begin{minipage}{120mm}
{\small \noindent{\bf Abstract.} For minimizers in a geometrically nonlinear
Cosserat model for micropolar elasticity of continua, we prove interior H\"older
regularity, up to isolated singular points that may be possible if the exponent
$p$ from the model is $2$ or in $(\frac{32}{15},3)$. 
The obstacle to full continuity turns out to be
the existence of certain minimizing homogeneous $p$-harmonic maps to $S^3$.
For those, we slightly improve existing regularity theorems in order to
achieve our result on the Cosserat model.\\\strut\\
{\bf MSC 2020.} 58E20; 74G40; 74B20.}
\end{minipage}
\end{center}

\section{Introduction and statement of results}

Cosserat micropolar elasticity is a framework for theories of continua
(as well as shells and rods) with some internal structure. The foundations
have been laid out by the brothers Eug\`ene and Fran\c{c}ois Cosserat
in 1909 \cite{CC}. In addition to stresses responding to translational
degrees of freedom of an elastic body, the framework also allows stresses
coming from rotational degrees of freedom assigned to every point
of the material.

On the other hand, $p$-harmonic maps are a well-established branch of
geometric analysis. They are critical points of the $p$-Dirichlet integral
$E^p(u):=\int_\Omega|Du|^p\,dx$ among mappings $u:\Omega\to N$, where
$N$ is some fixed Riemannian manifold.

The current paper exploits relations between both theories. Since the
$p$-Dirichlet integral, here for mappings to $SO(3)$, also appears in the
energy functional for useful models within Cosserat theory, the equations
for the latter couple the $p$-harmonic map equation with another one.
The analytic difficulties of Cosserat theories have their origin, at least
partially, in the geometric restriction of the rotational degrees of freedom
to $SO(3)$. But this is exactly the kind of restriction that has to be
understood in $p$-harmonic map theory, which has been developed to quite
some extent. We therefore aim at understanding the nonlinear aspects of
Cosserat theory better by using methods that have been successfully
established for $p$-harmonic maps. We address the question of
{\em partial regularity of minimizing weak solutions\/}, and find out
that methods invented by Luckhaus \cite{Lu} (based on \cite{SU1})
are particularly useful.

Many variants of Cosserat theory are available, and the author will not even
try to give an overview about the different models and the vast body of
results. Instead, we restrict to a particular instance of that theory
for micropolar elastic bodies, that has been studied in the framework
of the calculus of variations by Neff \cite{Ne1} and others. The elastic
body exists over a reference configuration that can be thought of a
subset $\Omega$ of $\R^3$. From that configuration, the body can be deformed,
shifting every point $x\in\Omega$ to some point $\p(x)\in\R^3$, such
that $\p(x)-x$ can be thought of its usually small dislocation.
Additionally we assume some structure of the material that attaches to
every $x\in\Omega$ an orthonormal frame that is free to rotate in $\R^3$
by an orthogonal matrix $R(x)\in SO(3)$. Both translations and rotations
cause material stresses, which are given by $R^tD\p(x)-I$ and $R^tDR$,
respectively. Here $I\in\R^{3\times3}$ is the identity matrix, and
we denote the transposed of the matrix $R$ by $R^t$. Note that $R^tDR$
is a $3$-tensor rather than a $2$-tensor, but since $R^t\d_iR\in so(3)$
for $i\in\{1,2,3\}$, there are only $9$ independent components. We will
not bother about aspects of modelling, and refer to the discussions in
\cite{Ne1} and \cite{NBO} instead.

Now let us describe the energy functional summing up the energy stored in
our elastic body. The contribution of the translation should be measuring
$R^tD\p-I$ somehow. The usual choice is
\[
  \mue\|\dev\sym(R^tD\p-I)\|_{L^2(\Omega)}^2
  +\muc\|\skew(R^tD\p-I)\|_{L^2(\Omega)}^2
  +3\muz\|\tr(R^tD-I)\|_{L^2(\Omega)}^2.
\]
Here $\dev\sym A$ is the deviatoric symmetric part $\frac12(A+A^t)-(\tr A)I$
of $A$, and $\skew A:=\frac12(A-A^t)$ is the skew-symmetric part.
Defining $P:\R^{3\times3}\to\R^{3\times3}$ to be the linear operator given by
\[
  PA:=\sqrt{\mue}\,\dev\sym A+\sqrt{\muc}\,\skew A
  +\sqrt{\muz}\,(\tr A)\,I,
\]
the term is written more simply as
\[
  \|P(R^tD\p-I)\|_{L^2(\Omega)}^2.
\]
The constants $\mue$, $\muc$, and $\muz$ will be assumed to be $>0$.
While this is completely usual for $\mue$ and $\muz$, it would be
desirable to allow the so-called ``Cosserat constant'' $\mu_c$ to
be $0$. Remember that elasticity theory usually involves $\sym D\p$
instead of $D\p$. This would, however, combine the ``geometric'' difficulty
for $R\in SO(3)$ with the ``coercivity issue'' for $\p$, and we are currently
not able to handle both. Anyway, the existence of minimizers has been
established in the $\mu_c>0$ case in \cite{Ne2} (where Cosserat is a special
case, discussed more explicitly in \cite{Ne1}). For $\mu_c=0$, on the
other hand, interesting
cases are open, see the discussion in \cite{Ne1}. And since we prefer dealing
with the regularity of minimizers in cases where they are known to exist, 
we have another reason to restrict to
$\mu_c>0$ in this paper. We expect that the generalization to $\mu_c\ge0$
would provide us with some interesting additional problems.

The contribution of the rotational stresses to the energy is simply
\[
  \lambda\|R^tDR\|_{L^p(\Omega)}^p
\]
for $\lambda>0$ and some parameter $p\ge2$.
On first glimpse, $p=2$ seems to be the
natural choice, but there are problems with the decoupling of linearized
equations that suggest that $p>2$ might be better for many purposes.
We shall see that even our regularity theory works slightly better
for $p$ larger than (but close to) $2$.
Some more general terms, involving parameters like $\mue,\muc,\muz$
above, have been proposed, but they seem less natural, since $R^t\d_iR$
is always skew-symmetric. (One can use $\Curl R$ instead of $DR$, though.)
Anyway, most of our regularity theory would work for those more general
energies, too, with the exception of those parts where point singularities
are removed. We therefore restrict to the simple term above. Since
$|R^tDR|=|DR|$, and since we can make one of the constants to be $1$,
we can even work with the simpler term
\[
  \|DR\|_{L^p(\Omega)}^p
\]
here, and of course this is the $p$-Dirichlet integral.

Our elastic body may be subject to exterior forces. Some of them,
e.g.\ mechanical ones, will act on the boundary of the body, only.
We need not consider them, because we are only concerned with
{\em interior regularity\/} in this paper. Some boundary regularity
may well be within reach of the methods presented here, but natural questions
seem to be more involved, like what happens at the edge between regions
with Dirichlet and Neumann boundary conditions.
Other forces, like gravity or electromagnetic forces, will act on
points in $\Omega$, and we have to account for such forces in our
functional. Exterior forces, given by a function $f:\Omega\to\R^3$,
are accounted for in the term
\[
  \int_\Omega(\p-x)\cdot f\,dx,
\]
while there may be also moments of force
affecting the rotational degrees of freedom,
given by $M:\Omega\to\R^{3\times3}$. They contribute to the energy via
\[
  \int_\Omega R\cdot M\,dx.
\]
For our domain $\Omega$ we have to assume that it is bounded. In order
to have existence of minimizers, we also assume that it is Lipschitz.

Summarizing, for a pair of functions $\p:\Omega\to\R^3$ and
$R:\Omega\to SO(3)$, we have the energy functional
\[
J(\p,R):=\int_\Omega\Big(|P(R^tD\p-I)|^2+|DR|^p
+(\p-x)\cdot f+R\cdot M\Big)\,dx.
\]
The topic of our paper is interior regularity of minimizers, which have
been proven to exist by Neff, given suitable boundary conditions.
Restricting to minimizers means we are only considering a static
problem here, no dynamics.
Minimizers are weak solutions of the Euler-Lagrange
equations for $J$, which are standard to derive. They read
\ea{elphi}
\div(RP^2(R^tD\p-I))&=&f,\\
\la{elR}
\Big(\div(|DR|^{p-2}DR)-\frac{2}{p}\,D\p\,P^2((D\p)^tR-I))-\frac1p\,M\Big)(x)
&\perp&T_{R(x)}SO(3).
\eea
The second one is an orthogonality relation, since variations of $R$
can only be made in directions tangential to $SO(3)$. Therefore,
\re{elR} represents only three independent equations rather than the
expected nine for the components of $R$. The ``missing'' six equations are
simply the requirement $R(x)\in SO(3)$.

The terms involving $f$ and $M$ are of lower order and different scaling
than the others, thus inflicting only minor complications to our study.
For all of this paper except this introduction and the last section,
we can therefore reasonably assume $f\equiv0$ and
$M\equiv0$. The minor changes necessary for nonvanishing $f$ and $M$
will be hinted at in Section~\ref{sect:forces}. 
For most of the paper, our functional therefore is
\[
  J(\p,R):=\int_\Omega\Big(|P(R^tD\p-I)|^2+|DR|^p\Big)\,dx.
\]

If we consider regularity for minimizers of $J$, we find that it cannot
be any better than the regularity of minimizing $p$-harmonic maps
$\Omega\to SO(3)$, i.e.\ minimizers of $E^p(R):=\int_\Omega|DR|^p\,dx$.
For minimizing $p$-harmonic maps of an $n$-dimensional domain $\Omega$,
partial regularity has been proven independently by Hardt/Lin \cite{HL}
and Fuchs \cite{Fu}, and with a more flexible proof by Luckhaus \cite{Lu}.
The result is that they are H\"older continuous (even $C^{1,\mu}$) in the
interior of $\Omega$ away from a closed set $\Sing(R)$ of Hausdorff dimension
$\le n-[p]-1$, and that $\Sing(R)$ is even discrete if $n-[p]-1=0$,
and empty if that is $<0$. For $n=3$, that means that one has
a discrete singular set for $p\in[2,3)$, while for $p\ge3$ we have
full H\"older continuity. We will find the same for minimizers
$(R,\p)$ of $J$, following Luckhaus' techniques and modifying them for
a functional where $\p$ and $R$ have different homogeneities and $\p$
is possibly unbounded, both of which are not allowed in \cite{Lu}.

Depending on the target, the singular set of $p$-harmonic maps
may be even smaller, or empty
altogether. In the $p=2$ case, there are results for special targets,
e.g.\ \cite{SU1} for minimizing harmonic maps to spheres. For $p>2$, 
there are some results, too, like \cite{XY} and \cite{Na}.

In our study, $R$ does not map to a sphere, but
to $SO(3)$, which however has $S^3$ as its universal cover.
This enables us to prove the following. Whenever $p>2$ and
there is a point singularity in some minimizer $(\p,R)$ for our
functional $J$, then there is also a minimizing $p$-harmonic
$u:B^3\to S^3$ having a point singularity. Unfortunately, the results in
\cite{XY} or \cite{Na} are not strong enough to exclude the latter.
But we can do so, at least for $p\in(2,\frac{32}{15}]$. 

To do this, we
partially follow recent progress by Chang, Chen, and Wei \cite{CCW}
for $p$-harmonic functions
to $\R$, resulting in an ``improved Kato inequality''. Our results
strongly depend on the values of constants in estimates and therefore
are almost certainly far from optimal. But anyway, this shows that
for $p\in(2,\frac{32}{15}]$ (and for $p=3$, and almost trivially for $p>3$)
minimizers of our Cosserat energy are H\"older continuous on the interior
of the domain. This leaves the possibility of point singularities
only for $p=2$ or $p\in(\frac{32}{15},3)$.

Now, here are our precise results.
The following theorem combines the statements of the Propositions \ref{hoel1},
\ref{hoel2}, \ref{reg1}, \ref{reg2}, and \ref{reg3} formulated and proven
below in this paper.

\begin{theorem}[interior (partial) regularity for minimizers]
  Assume $\mue,\muc,\muz>0$, and $p\ge2$. Let $\Omega\subset\R^3$
  be a bounded open domain, let functions
  $f\in C^{0,\mu}(\overline{\Omega},\R^3)$ and
  $M\in C^0(\overline{\Omega},\R^{3\times3})$ be prescribed, and let
  $(\p,R)\in W^{1,2}(\Omega,\R^3)\times W^{1,p}(\Omega,SO(3))$ be a minimizer
  of the functional $J$. Then there is a discrete subset $\Sing(\p,R)$
  of $\Omega$ such that
  \[
  (\p,R)\in C^{1,\mu}_\loc(\Omega\setminus\Sing(\p,R),\R^3)
  \times C^{0,\mu}_\loc(\Omega\setminus\Sing(\p,R),SO(3))
  \]
  for every $\mu\in(0,\frac2p)$ if $f\equiv0$, and for every
  $\mu\in(0,\frac1{2p})$ if $f\not\equiv0$.

  Moreover, $\Sing(\p,R)$ is empty, and therefore $\p$, $D\p$, and $R$
  locally H\"older continuous on all of $\Omega$, if one of the following
  conditions holds.

  (i) $p=2$ and $\mue=\muc=\muz$,

  (ii) $p\in(2,\frac{32}{15}]$,

  (iii) $p\ge3$.
\end{theorem}

The obvious question this theorem provokes is if there can be singular points
at all for minimizers of $J$. We cannot answer this question exactly, but,
for every $p\in[2,3)$, we do find an explicit weak solution for the system
of the Euler-Lagrange equations that has a singular point 
--- one more motivation to study regularity
theory. We do not know if our example is minimizing for any $p$.

Our observations are reflected by the regularity theory for
$p$-harmonic maps, where, too, much more is known about minimizers
than about more general weak solutions. In dimensions $>2$, one cannot
expect any good regularity in general, since Rivi\`ere \cite{Ri} has
constructed a weakly harmonic map $B^3\to S^2$ that is discontinuous
in every point of $B^3$. Under the additional assumption of
stationarity, weakly $p$-harmonic maps are sometimes known to be
H\"older continuous outside a closed set of vanishing
$(n-p)$-dimensional Hausdorff measure, proven by Bethuel \cite{Be} for
$p=2$ and any compact target, and by Toro and Wang \cite{TW} for
general $p$, but only if the target is a compact homogeneous
space. Here, a weak solution is called stationary if it is also
critical with respect to variations in the domain.

Assuming that our Cosserat system
has regularity just as good as for harmonic maps, we would expect
a singular set of vanishing $(3-p)$-dimensional Hausdorff measure for
stationary weak solutions, and no good regularity theory at all for
just weak solutions. In this paper, however, we only consider minimizers.

The paper is organized as follows. In Section~\ref{sect:ex}, we demonstrate
by our singular weak solution that regularity theory is an issue at all.
In Section~\ref{sect:three}, we show that in the simple case $p=2$ and
$\mue=\muc=\muz$ Luckhaus' result from \cite{Lu} already gives some
partial regularity for the Cosserat body, but not as much regularity
as we will achieve here. In Section~\ref{sect:four}, we adapt Luckhaus'
proof of partial regularity by blowing up in the target and domain in order
to compare with a simplified system. This results in an $\eps_0$-regularity
theorem stating that singularities can only occur in points where enough
energy concentrates. Section~\ref{sect:five}, also inspired by the
techniques for ($p$-)harmonic maps, features a monotonicity formula saying
that minimizers are automatically in some Morrey space rather than just
in $W^{1,2}\times W^{1,p}$. We apply this in order to get regularity up to
the possibility of isolated singularities. The obstruction to full regularity
is identified to be the existence of certain $p$-harmonic maps to $S^3$.
Therefore, in Section~\ref{sect:six}, we try to exclude their existence under
suitable assumptions, which leads to full H\"older regularity for some 
exponents $p$.
Finally, Section~\ref{sect:forces} discusses the changes necessary to allow
for exterior forces and moments, that in the sections before were
assumed to vanish.

{\bf Acknowledgement. }This paper was written while the author was visiting
the group of Prof.~Frank Duzaar at Erlangen. The author would like to thank
him and his group for hospitality and support.

\section{A singular example}\la{sect:ex}

From now on, $f\equiv0$ and $M\equiv0$ will be assumed until we discuss the
case of their nonvanishing in Section~\ref{sect:forces}.

If all constants agree, by which we mean $\mue=\muc=\muz=1$,
the Euler-Lagrange equations simplify according to $P=\id$, giving
\eas
\Delta\p-\div R&=&0,\\
\div(|DR|^{p-2}DR)+\frac{2}{p}\,D\p&\perp&T_RSO(3),
\eeas
where in the latter we have removed the term $D\p\,D\p^t\,R$ since it is always
orthogonal to $T_RSO(3)$.
In this simple case, we can write down an explicit weak solution that
exhibits a singularity. It is given by $(\p,R):B^3\to\R^3\times SO(3)$,
\eas
\p(x)&:=&\frac43\,x\log|x|,\\
R(x)&:=&\frac2{|x|^2}\,x\otimes x-I=\frac1{|x|^2}\left(\begin{array}{ccc}
  x_1^2-x_2^2-x_3^2 & 2x_1x_2 & 2x_1x_3 \\
  2x_1x_2 & x_2^2-x_1^2-x_3^2 & 2x_2x_3 \\
  2x_1x_3 & 2x_2x_3 & x_3^2-x_1^2-x_2^2
\end{array}\right)\,.
\eeas
We perform a few calculations to see that this is a solution on
$B^3\setminus\{0\}$. We have
\eas
\d_i\p(x)&=&\frac43\,\Big(e_i\log|x|+\frac{x_ix}{|x|^2}\Big),\\
\d_i^2\p(x)&=&\frac43\,\Big(2\,\frac{x_ie_i}{|x|^2}+\frac{x}{|x|^2}
-2\,\frac{x_i^2x}{|x|^4}\Big),\\
\Delta\p(x)&=&4\,\frac{x}{|x|^2}\,,
\eeas
and
\eas
\div R&=&4\,\frac{x}{|x|^2}\,,\\
\d_kR_{ij}&=&\frac2{|x|^2}\,(\del_{jk}x_i+\del_{ik}x_j)
  -\frac{4x_ix_jx_k}{|x|^4}\,,\\
|DR|^2&=&\frac4{|x|^2}\,,\\
\div(|DR|^{p-2}DR)&=&\frac{2^p}{|x|^p}\,
  \Big(I-\frac3{|x|^2}\,x\otimes x\Big).
\eeas
Now we immediately read off the first Euler-Lagrange equation. And for
every $x\in\R^3$, we see that $\div(|DR(x)|^{p-2}DR(x))$ and $D\p(x)$ are
linear combinations of $I$ and $x\otimes x$. A matrix $A\in\R^{3\times3}$
is perpendicular to $T_RSO(3)$ if $R^tA$ is perpendicular to $so(3)$,
which holds if and only if $R^tA$ is symmetric. But $R^t(x)I=R^t(x)$ and
$R^t(x)(x\otimes x)=x\otimes x$ are both symmetric, which means that
$\div(|DR|^{p-2}DR)+\frac{2}{p}\,D\p$ is perpendicular to $T_{R(x)}SO(3)$
for every $x\in B^3\setminus\{0\}$. Therefore also the second equation holds
away from the origin.

Note that $(\p,R)\in W^{1,2}\times W^{1,p}$ if
$2\le p<3$, and we easily find that, for those $p$, $(\p,R)$ is a weak
solution of the Euler-Lagrange equations. Of course, it is smooth on
$B^3\setminus\{0\}$. In $0$, $R$ is not even continuous, while
$\p$ is H\"older continuous for every exponent $<1$, but not differentiable.

The example shows that we must be ready to expect at least point singularities
for weak solutions of our model, as long as $p<3$. The solutions constructed
by Neff, on the other hand, are better than just weak solutions, they are
minimizers of $J$. As often in the calculus of variations, we will find
that sometimes minimizers have better regularity properties than other
weak solutions.

\section{A quick application of a result by Luckhaus}\la{sect:three}

We first consider the case $p=2$ and show that it is within the framework
of a paper \cite{Lu} by Luckhaus which was written with focus on
$p$-harmonic maps.

Assume that $(\p,R)\in W^{1,2}(\Omega,\R^3\times SO(3))$ is a minimizer
of $J$ (subject to suitable boundary conditions). We check that assumptions
from \cite{Lu} are fulfilled. To this end, we introduce some notation.
We write $N:=\R^3\times SO(3)$, $x$ for the independent
variable in $\Omega$ and $y=(y_1,y_2)\in N$ for the
dependent variable holding $(\p,R)$. Moreover, $z=(z_1,z_2)$ stands for
the variable that $(D\p,DR)$ take their values in. This means
\[
J(\p,R)=\int W((\p(x),R(x)),(D\p(x),DR(x)))\,dx,
\]
where here
\[
W(y,z):=\mue|\dev\sym y_2^tz_1|^2+\muc|\skew y_2^tz_1|^2+\muz
\,|\tr y_2^tz_1-3|^2+|z_2|^2.
\]
Since $W$ does not depend explicitly on $x$, the situation is even slightly
simpler than Luckhaus'.

Luckhaus has the following conditions (adapted here for $x$-independent
functionals)
\begin{eqnarray*}
  (A1a)&&c^{-1}|z|^2-1\le W(y,z)\le c|z|^2+1,\\
  (A1b)&&\lim_{y\to y_0}\sup_z(1+|z|)^{-2}|W(y,z)-W(y_0,z)|=0,\\
  (A1c)&&W(y,z)\mbox{ is convex in $z$ for all $y$.}
\end{eqnarray*}
Defining for $y\in N$ the set
\[
  {\mathcal H}_y:=\{F:\,\mbox{there exist $\a_i\to\infty$ s.t.\ }
      F(z)=\lim\a_i^{-2}W(y,\a_iz)\}
\]
as well as
\[
  {\mathcal H}_\infty:=\{F:\,\mbox{there exist $\a_i\to\infty$, $|y_i|\to\infty$
    s.t.\ }F(z)=\lim\a_i^{-2}W(y_i,\a_iz)\},
\]
Luckhaus has one more condition,  
\begin{quote}
  $(A2)$ If $F\in{\mathcal H}_y$ or $F\in{\mathcal H}_\infty$ and
  $v\in W^{1,2}(B_1,T)$ for $T:=T_yN$ or $T:=\lim T_{y_i}N$, and
  $$\div(D_zF(Dv))=0,$$
  then $v$ is $\mu$-H\"older continuous for some $\mu\in (0,1]$.
\end{quote}
For the target manifold $N$, in \cite{Lu} it is enough for part (a) of his
theorem, that the nearest point retraction when restricted to an
$\epsilon$-neighborhood of $N$ has a Lipschitz constant approaching $1$
uniformly as $\eps\searrow0$. This is clearly fulfilled for our
$N=\R^3\times SO(3)$ because of its bounded curvatures.

Hence it remains to check the $(A1)$ and $(A2)$ assumptions.
We see that $(A1a)$ clearly holds since for $c:=\max\{\muc,\mue,\muz,1\}$
we have
\[
c^{-1}|z|^2-18\le W(y,z)\le c|z|^2+18,
\]
and $18$ is certainly as good as $1$ for our purposes.
The condition $(A1b)$ is immediately read off for our functional,
and the convexity required for $(A1c)$ also clearly holds. 

We now observe that ${\mathcal H}_y$ consists of the single function
\[
F(z):=\mue|\dev\sym y_2^tz_1|^2+\muc|\skew y_2^tz_1|^2+\muz
\,|\tr y_2^tz_1|^2+|z_2|^2.
\]
And ${\mathcal H}_\infty$ consists of all those for all $y_2\in SO(3)$.
Every such function is a homogeneous quadratic form in $z$ which is
positive definite. Therefore, $\div D_zF(Dv)=0$ is an elliptic equation for
$v$ with constant coefficients, whose solutions are of course H\"older
continuous. This shows that $(A2)$ also holds.

We can therefore apply Theorem (a) from \cite{Lu} and conclude that
any minimizer $(\p,R)$ of $J$ is H\"older continuous on the interior of
$\Omega$ outside a closed singular set $\Sigma$ for which we
have $H^1(\Sigma)=0$, where $H^1$ means the $1$-dimensional Hausdorff
measure.

A second theorem from \cite{Lu} states that the singular set consists
of isolated points only, but it cannot be applied immediately, because it
requires a compact target manifold $N$, plus more assumptions. Our
$N=\R^3\times SO(3)$ is not compact, but the noncompact factor is
quite trivial, which allows us to work around this in the sections that
follow.

\section{A first partial regularity result for general $p\ge2$}
\la{sect:four}

If $p>2$, the reasoning from the previous section does not work
properly, and we have to use modifications of Luckhaus'
arguments from \cite{Lu} to prove partial H\"older continuity.
Since the unknown functions $\p\in W^{1,2}(\Omega,\R^3)$ and
$R\in W^{1,p}(\Omega,\R^n)$ now are in Sobolev spaces of different
scaling, \cite{Lu} cannot be applied immediately. The overall strategy,
however, continues to work, and we will be able to use the key ideas,
including the ``Luckhaus Lemma'', without too much modification.
Our first lemma is a discrete version of Morrey's Dirichlet growth
condition as a first step to partial regularity. We globally assume the
constants $\mue,\muc,\muz>0$ to be fixed, in fact most constants will
depend on those.

Morrey's Dirichlet growth criterion would imply local H\"older continuity
of $R$ once we know that $\rho^{1-3/p}\|DR\|_{L^p(B_\rho(x_0))}$ is bounded
by some $c\rho^\mu$ for all $B_\rho(x_0)\subset\Omega$. The following
Lemma establishes a discrete version of this which allows the same conclusion.

\begin{lemma}[discrete Morrey condition]\la{morr}
  We fix $\mu\in(0,1)$. Then
  there exists constants $\eps_0>0$, $\th\in(0,1)$, and $\rho_0\in(0,1)$
  such that the following holds for every minimizer
  $(\p,R)\in W^{1,2}(B^3,\R^3)\times W^{1,p}(B^3,SO(3))$
  of $J$ subject to its boundary data. For every ball $B_{\rho_0}(x_0)
  \subset B^3$ with $|x_0|\le\frac12$ and any $\rho\in(0,\rho_0)$,
  the condition
  \eas
  \rho^{2\mu}\le
  \rho^{-1}\|D\p\|_{L^2(B_\rho(x_0))}^2+\rho^{p-3}\|DR\|_{L^p(B_\rho(x_0))}^p
  \le\eps_0
  \eeas
  imply
  \eas
  \lll{(\th\rho)^{-1}\|D\p\|_{L^2(B_{\th\rho}(x_0))}^2
    +(\th\rho)^{p-3}\|DR\|_{L^p(B_{\th\rho}(x_0))}^p}\\
  &\le&\th^{2\mu}(\|D\p\|_{L^2(B_\rho(x_0))}^2+\rho^{p-3}\|DR\|_{L^p(B_\rho(x_0))}^p).
  \eeas
\end{lemma}

{\bf Proof.} Assume that the assertion does not hold. Then there are balls
$B_{\rho_i}(x_i)$ with $|x_i|\le\frac12$ and $\rho_i\searrow0$ such that
\eas
\g_i&:=&\rho_i^{-1/2}\|D\p\|_{L^2(B_{\rho_i}(x_i))}\searrow0,\\
\del_i&:=&\rho_i^{1-3/p}\|DR\|_{L^p(B_{\rho_i}(x_i))}\searrow0.
\eeas
and $\rho_i^{2\mu}\le\g_i^2+\eps^p$, but 
\eas
\lll{(\th\rho_i)^{-1}\|D\p\|_{L^2(B_{\th\rho_i}(x_i))}^2
  +(\th\rho_i)^{p-3}\|DR\|_{L^p(B_{\th\rho_i}(x_i))}^p}\\
&>&\th^{2\mu}(\|D\p\|_{L^2(B_{\rho_i}(x_i))}^2
+\rho_i^{p-3}\|DR\|_{L^p(B_{\rho_i}(x_i))}^p).
\eeas

Now we do a suitable rescaling and define $(\p_i,R_i)\in W^{1,2}(B^3,\R^3)
\times W^{1,p}(B^3,N_i)$ by
\eas
\p_i(x)&:=&\g_i^{-1}(\p(x_i+\rho_i x)-\opi),\\
R_i(x)&:=&\del_i^{-1}(R(x_i+\rho_i x)-\oRi),
\eeas
where here $\opi$ and $\oRi$ are the mean values
\[
\opi:=\mint_{B_{\rho_i}(x_i)}\p\,dx,\qquad
\oRi:=\mint_{B_{\rho_i}(x_i)}R\,dx,
\]
and $N_i$ is the rescaled and shifted target manifold $N_i:=\del_i^{-1}(N-\oRi)$.
Then $(\p_i,R_i)$ minimizes the rescaled functional
\[
J_i(\tp,\tR):=\frac{\g_i^2}{\g_i^2+\del_i^p}
\int_{B^3}|P((\oRi+\del_i\tR)^tD\tp-\rho_i\g_i^{-1}I)|^2\,dx
+\frac{\del_i^p}{\g_i^2+\del_i^p}
\int_{B^3}|D\tR|^p\,dx.
\]
Here the denominator is chosen such that,
after passing to a subsequence, we may assume
$\frac{\g_i^2}{\g_i^2+\del_i^p}\to\sigma$ and
$\frac{\del_i^p}{\g_i^2+\del_i^p}\to1-\sigma$
for some $\sigma\in[0,1]$. Note also that
$\frac{\g_i^2}{\g_i^2+\del_i^p}(\rho_i\g_i^{-1})^2\to0$ since
$\rho_i^{2\mu}/(\g_i^2+\del_i^p)\to 0$ and $\mu<1$.
And the sequence $\oRi$ is certainly bounded, since
it cannot leave the convex hull of $SO(3)$ in $\R^{3\times3}$. We even have
\[
\dist(\oRi,SO(3))\le c\rho_i^{-3/p}\|R-\oRi\|_{L^p(B_{\rho_i}(x_i))}
\le c\rho_i^{1-3/p}\|DR\|_{L^p(B_{\rho_i}(x_i))}=c\del_i\to0,
\]
and therefore we can assume $\oRi\to T$ for some $T\in SO(3)$.

We now consider the corresponding subsequence of our $J_i$-minimizers
$(\p_i,R_i)$ and hope that it converges (in a suitable sense)
to a minimizer of the limit functional
\[
J_\infty(\tp,\tR):=\sigma\int_{B^{3}}|P(T^tD\tp)|^2\,dx
+(1-\sigma)\int_{B^3}|D\tR|^p\,dx.
\]
This may be not quite true, but almost, since it holds away from $\d B^3$.
This is stated and proved in Lemma \ref{comp} below. From that we know
that on every compact subset of $B^3$, $(\p_i,R_i)$ converges in
$W^{1,2}\times W^{1,p}$-norm to a minimizer $(\p_\infty,R_\infty)$ of $J_\infty$.
Note that the $N_i$ converge locally in Hausdorff distance
to some $3$-dimensional subspace of $\R^{3\times3}$ which we denote by $N_\infty$.
Actually, $N_\infty$ is the limit of the affine tangent spaces
$\oRi+T_{\oRi}N_i$, and $0\in N_\infty$ because all $R_i$ have mean value $0$.

In preparation of Lemma \ref{comp}, we note that we have
weak (sub-)convergence $(\p_i,R_i)\rightharpoonup(\p_\infty,R_\infty)$.
This holds because $\|D\p_i\|_{L^2(B^3)}=1$ and $\|DR_i\|_{L^p(B^3)}=1$ by
construction, and both $\p_i$ and $R_i$ have mean values $0$.

Up to passing to another subsequence, we then have
\eq{co}
\oRi+\del_iR_i\to T
\eeq
pointwise almost everywhere, for the constant
matrix $T$ found above. We see this by combining the pointwise convergences
$\oRi\to T$, $R_i\to R_\infty$ and $\del_i\to0$ in
\[
\oRi+\del_i R_i=\oRi+\del_i R+\del_i(R_i-R)\to\lim_{i\to\infty}\oRi=T.
\]
This explains why we consider $J_\infty$ the right limit functional.

Having applied Lemma \ref{comp}, we return to our proof of Lemma \ref{morr},
and we have to distinguish three cases. In what follows, $J_{i,r}$ means
the same functional as $J_i$, but with integration over $B_r$ instead
of $B^3$.

The first case is $\sigma\in(0,1)$. Then $(\p_\infty,R_\infty)$ can only minimize
$J_{\infty,r}$ if $\p_\infty$ is a minimizer of $\int_{B_r}|P(T^tD\tp)|^2\,dx$ and
$R_\infty$ is a minimizer of $\int_{B_r}|D\tR|^p\,dx$. Then $\p_\infty$ solves an
elliptic system with constant coefficients, and hence is H\"older continuous.
And $R_\infty$ (now with values in a vector space, contrary to $R$)
solves the $p$-Laplace system considered by Uhlenbeck \cite{Uh} who also
proved H\"older continuity. More precisely, we have the regularity estimates
\eas
\int_{B_\th}|D\p_\infty|^2\,dx&\le&c\th^{2\nu}\th\int_{B_1}|D\p_\infty|^2\,dx,\\
\int_{B_\th}|DR_\infty|^p\,dx&\le&c\th^{p\nu}\th^{3-p}\int_{B_1}|DR_\infty|^p\,dx
\eeas
for all $\nu\in(0,1)$, $\th<\frac12$, with $c$ depending only on $\nu$
(and $P$, which is considered fixed). The first of the estimates is a
well-known standard estimate. Note that the equation $\p_\infty$ solves
depends on $T$, but the constant $c$ does not, since it only depends on
the ellipticity constant of the operator, i.e.\ only on $P$.
The second estimate follows from Uhlenbeck's H\"older regularity result
(saying that the $p$-harmonic $R_\infty$ is in $C^{0,\mu}$ for every
$\mu\in(0,1)$, since it is
even in $C^{1,\a}$), using \cite[Lemma 2(a)]{Lu}.

We now choose $\th$ small enough in
order to reach the desired contradiction. By norm convergence on $B_\th$
and weak convergence on $B_1=B^3$, we find, for $\nu:=\frac{\mu+1}2$,
\ea{it}
\lll{\lim_{i\to\infty}\frac1{\g_i^2+\del_i^p}\,
  ((\th\rho_i)^{-1}\|D\p\|_{L^2(B_{\th\rho_i}(x_i))}^2
  +(\th\rho_i)^{p-3}\|DR\|_{L^p(B_{\th\rho_i}(x_i))}^p)}\nn\\
&=&\lim_{i\to\infty}\frac1{\g_i^2+\del_i^p}\,(\th^{-1}\g_i^2\|D\p_i\|_{L^2(B_\th)}^2
+\th^{p-3}\del_i^p\|DR_i\|_{L^p(B_\th)}^p)\nn\\
&=&\sigma\th^{-1}\|D\p_\infty\|_{L^2(B_\th)}^2
+(1-\sigma)\th^{p-3}\|DR_\infty\|_{L^p(B_\th)}^p\nn\\
&\le&c\th^{2\nu}(\sigma\|D\p_\infty\|_{L^2(B_1)}^2
+(1-\sigma)\|DR_\infty\|_{L^p(B_1)}^p)\nn\\
&\le&c\th^{2\nu}\nn\\
&\le&\th^{2\mu}
\eea
if $\th$ is taken small enough to have $c\th^{1-\mu}\le1$. This contradicts
our original assumption in the first case.

The second case is $\sigma=0$, in which case we have a minimizer
$R_\infty$ of $\int_{B_r}|D\tR|^p\,dx$, but the information on
$\p_\infty$ has been lost in the limit. However, we still have the estimate
\re{it}, since the $\|D\p_\infty\|$-term we now cannot estimate has the
coefficient $0$, anyway. The third case, $\sigma=1$, uses the same arguments
with $R$ taking the role of $\p$.

This proves the Lemma, up to Lemma \ref{comp} below.\qed

The compactness lemma is proven almost exactly as in Luckhaus' paper,
we formulate a sketchy proof only for the reader's convenience. 

\begin{lemma}[compactness]\la{comp}
  In the proof of Lemma \ref{morr}, the weakly convergent sequence
  $(\p_i,R_i)$ of $J_i$-minimizers converges even in
  $W^{1,2}\times W^{1,p}$ on every compact subset of $B^3$, and the
  limit $(\p_\infty,R_\infty)$ minimizes the limit functional $J_\infty$.
\end{lemma}

{\bf Sketch of proof.}
Let $\psi\in W^{1,2}(B^3,\R^3)$ and $Q\in W^{1,p}(B^3,N_\infty)$ be given
such that $\psi=\p_\infty$ and $Q=R_\infty$ on some neighborhood of
$\d B^3$. Then (...) there exists $r$ close to $1$
such that $\psi=\p_\infty$ and $Q=R_\infty$ almost everywhere on
$\d B_r$, and there exist $Q_i\in W^{1,p}(B_r,N_i)$ such that
$Q_i\to Q$ in $W^{1,p}(B^3,\R^{3\times3})$,
\eq{c1}
\|D\psi\|_{L^2(\d B_r)}+\sup_i(\|D\p_i\|_{L^2(\d B_r)}+\|DQ_i\|_{L^p(\d B_r)}
+\|DR_i\|_{L^p(\d B_r)})=:K<\infty
\eeq
and
\eq{c2}
\lim_{i\to\infty}\|Q_i-R_i\|_{L^p(\d B_r)}=0.
\eeq
Now the ``Luckhaus Lemma'' (Lemma~1 from [Lu])
provides us with two sequences of functions
$\zeta_i\in W^{1,2}(B^3,\R^3)$ and $P_i\in W^{1,p}(B^3,\R^{3\times3})$
as well as number sequences $\lambda_i\searrow0$ and $r_i\searrow0$
such that
\ea{c3}
&&\zeta_i=\p_i\mbox{ and }P_i=R_i\mbox{ on }B^3\setminus B_r,\\
\la{c4}
&&\zeta_i(x)=\psi(\tf{x}{1-\lambda_i})\mbox{ and }P_i(x)
=Q_i(\tf{x}{1-\lambda_i})\mbox{ on }B_{(1-\lambda_i)r},\\
\la{c5}
&&P_i(B^3)\subseteq B_{r_i}(N_i),\\
\la{c6}
&&\|D\zeta_i\|_{L^2(B_r\setminus B_{(1-\lambda_i)r})}
+\|DP_i\|_{L^p(B_r\setminus B_{(1-\lambda_i)r})}\le c\l_i^{1/p}K.
\eea
These equations together with the assumptions give that
the $\zeta_i$ are bounded in $W^{1,2}(B^3,\R^3)$ and the $P_i$ are bounded
in $W^{1,p}(B^3,\R^{3\times3})$.

Denote by $J_{i,r}$ the variant of $J_i$ where integration is over
$B_r$ instead of $B^3$. We apply the pointwise convergence in \re{co}
together with the dominated convergence theorem in the first step, and
weak lower semicontinuity in the second to infer
\ea{c7}
J_{\infty,r}(\p_\infty,R_\infty)&=&\lim_{i\to\infty}
\frac{\g_i^2}{\g_i^2+\del_i^p}
\int_{B_r}|P((\oRi+\del_iR_i)^tD\p_\infty-\rho_i\g_i^{-1}I)|^2\,dx\nn\\
&&+\lim_{i\to\infty}\frac{\del_i^p}{\g_i^2+\del_i^p}
\int_{B_r}|DR_\infty|^p\,dx.\nn\\
&\le&\liminf_{i,j\to\infty}
\frac{\g_i^2}{\g_i^2+\del_i^p}
\int_{B_r}|P((\oRi+\del_iR_i)^tD\p_j-\rho_i\g_i^{-1}I)|^2\,dx\nn\\
&&+\liminf_{i,j\to\infty}\frac{\del_i^p}{\g_i^2+\del_i^p}
\int_{B_r}|DR_j|^p\,dx.\nn\\
&\le&\liminf_{i\to\infty}J_{i,r}(\p_i,R_i).
\eea
Now we use that $(\p_i,R_i)$ is a $J_i$-minimizer and has the same boundary
values on $\d B_r$ as $(\zeta_i,\pi_i\circ P_i)$, where here $\pi_i$ is the
nearest point retraction onto $N_i$. Since the $N_i$ are
magnifications of $SO(3)$, the Lipschitz constants of $\pi_i$ (restricted
to tubes of width $1$, say, around $N_i$) are bounded independently on $i$.
And $\pi_i\circ P_i$ differs from $P_i$ only on the annuli $B_r\setminus
B_{(1-\lambda_i)r}$ where the $p$-energy of $P_i$ approaches $0$. Hence
the $p$-energies of $\pi_i\circ P_i$ and $P_i$ differ only by $o(1)$,
and the same applies for the $p$-energies of $P_i$ and $Q_i$, as well
as the $2$-energies of $\psi$ (independent from $i$) and $\zeta_i$.
For example, we have
\eas
  \lll{\frac{\g_i^2}{\g_i^2+\del_i^p}\int_{B_r}\Big(|P((\oRi+\del_iR_i)^tD\p_i
    -\rho_i\g_i^{-1}I)|^2}\\
  \lll{\qquad\quad-|P((\oRi+\del_i\pi_i\circ P_i)^tD\p_i
    -\rho_i\g_i^{-1}I)|^2\Big)\,dx}\\
  &\le&c\,\frac{\g_i^2}{\g_i^2+\del_i^p}\int_{B_r\setminus B_{1-\l_ir}}
    (|D\p_i|^2+\rho_i^2\g_i^{-2})\,dx\\
  &=&c\,\frac{\rho_i^2}{\g_i^2+\del_i^p}\int_{B_{r\rho_i}(x_0)\setminus
    B_{(1-\l_i)r\rho_i}(x_0)}(|D\p|^2+1)\,dx\\
  &\le&c\int_{B_{r\rho_i}(x_0)\setminus B_{(1-\l_i)r\rho_i}(x_0)}(|D\p|^2+1)\,dx\\
  &\to&0.
\eeas  
Together with similar estimates, we find
$J_{i,r}(\p_i,R_i)=J_{i,r}(\zeta_i,\pi_i\circ P_i)+o(1)$. We combine this with
the minimality of $(\p_i,R_i)$ to continue estimate \re{c7} in
\ea{c8}
J_{\infty,r}(\p_\infty,R_\infty)
&\le&\liminf_{i\to\infty}J_{i,r}(\p_i,R_i)\nn\\
&\le&\liminf_{i\to\infty}J_{i,r}(\zeta_i,\pi_i\circ P_i)\nn\\
&\le&\liminf_{i\to\infty}J_{i,r}(\psi,Q_i)\nn\\
&=&J_{\infty,r}(\psi,Q),
\eea
where in the last step we have used $Q_i\to Q$ in $W^{1,p}$ and a similar
reasoning as in \re{co} which shows $\oRi+\del_iQ_i\to T$ pointwise almost
everywhere. This proves that $(\p_\infty,R_\infty)$ is $J_{\infty,r}$-minimizing
with respect to its boundary values on $\d B_r$. Since $r<1$ can be chosen
arbitrarily close to $1$, we have that $(\p_\infty,R_\infty)$ is
$J_{\infty,r}$-minimizing on every compact subset of $B^3$.

It is allowed to chose $(\psi,Q)=(\p_\infty,R_\infty)$ in \re{c8}, and this gives
\[
\lim_{i\to\infty}J_{i,r}(\p_i,R_i)=J_{\infty,r}(\p_\infty,R_\infty),
\]
which by strict convexity is easily seen to imply convergence in
$W^{1,2}\times W^{1,p}$ on $B_r$. This completes the proof of
Lemma~\ref{comp}.\qed

It is now relatively standard to proceed from Lemma \ref{morr} 
to the following partial H\"older regularity statement.

\begin{proposition}[$\eps_0$-regularity and partial H\"older
    continuity]\la{hoel1}
  Assume $2\le p\le3$, $\mu\in(0,\frac2p)$, and the assumptions made above.
  Then every $J$-minimizer $(\p,R)\in
  W^{1,2}(\Omega,\R^3)\times W^{1,p}(\Omega,SO(3))$ is locally in $C^{1,\mu}\times
  C^{0,\mu}$ on $\Omega\setminus\Sigma$, where
  \[
  \Sigma:=\{x_0\in\Omega:\liminf_{\rho\searrow0}
  \rho^{-1}\|D\p\|_{L^2(B_\rho(x_0))}^2+\rho^{3-p}\|DR\|_{L^p(B_\rho(x_0))}^p\ge\eps_0\}
  \]
  for some sufficiently small $\eps_0>0$.
  
  The set $\Sigma$ is relatively
  closed in $\Omega$, and $\HH^1(\Sigma)=0$.
\end{proposition}

{\bf Proof. }The H\"older continuity of $R$ on $\Omega\setminus\Sigma$
is more or less exactly Luckhaus' argument which is as follows.
For the moment, let $\mu\in(0,1)$.
For every $x_0\in\Omega\setminus\Sigma$, there is some $s>0$ such that
\eq{env}
s^{-1}\|D\p\|_{L^2(B_s(x_0))}^2+s^{p-3}\|DR\|_{L^p(B_s(x_0))}^p\le2\eps_0.
\eeq
The key is to prove the energy estimate
\eq{ene}
r^{-1}\|D\p\|_{L^2(B_r(x))}^2+r^{p-3}\|DR\|_{L^p(B_r(x))}^p\le C\eps_0\Big(\frac{r}s
\Big)^{2\mu}
\eeq
for every $x\in B_{s/2}(x_0)$ and every $r\in(0,\frac{s}2)$.
Once we have that, the H\"older continuity of $R$ on $B_{s/2}(x_0)$
follows using Morrey's Dirichlet growth criterion.

To prove \re{ene}, we use Lemma \ref{morr}, the $\eps_0$ of which we
denote by $\tilde{\eps}_0$ here. 
Let $s_0:=\min\{\frac{s}2,\th^{(p-3)/2\mu-1}\tilde{\eps}_0^{1/2\mu}\}$.
By \re{env}, we have
\[
s_0^{-1}\|D\p\|_{L^2(B_{s_0}(x))}^2+s_0^{p-3}\|DR\|_{L^p(B_{s_0}(x))}^p
\le C_0\eps_0=\tilde{\eps}_0
\]
for every $x\in B_{s/2}(x_0)$ if we have chosen $\eps_0$ accordingly.
We abbreviate $\Phi(r):=r^{-1}\|D\p\|_{L^2(B_r(x))}^2+r^{p-3}\|DR\|_{L^p(B_r(x))}^p$.
By induction, we prove the claim $\Phi(\th^ks_0)\le\th^{2k\mu}\tilde{\eps}_0$
for all $k\in\N$. For this clearly holds for $k=0$, and let us assume that
$\Phi(\th^{k-1}s_0)\le\th^{2(k-1)\mu}\tilde{\eps}_0$ has already been proved. Then
either $(\th^{k-1}s_0)^{2\mu}<\Phi(\th^{k-1}s_0)$, and Lemma \ref{morr}
implies
\[
\Phi(\th^ks_0)\le\th^{2\mu}\Phi(\th^{k-1}s_0)\le\th^{2k\mu}\tilde{\eps}_0,
\]
or $(\th^{k-1}s_0)^{2\mu}\ge\Phi(\th^{k-1}s_0)$, and hence
\[
\Phi(\th^ks_0)\le\th^{p-3}\Phi(\th^{k-1}s_0)\le\th^{p-3}(\th^{k-1}s_0)^{2\mu}
=\th^{2k\mu}\th^{p-3-2\mu}s_0^{2\mu}\le\th^{2k\mu}\tilde{\eps}_0.
\]
Note that the smallness condition of Lemma \ref{morr} is fulfilled in
every step. We have thus proven $\Phi(\th^ks_0)\le\th^{2k\mu}\tilde{\eps_0}$
for all $k$, and this clearly implies \re{ene}, and hence the asserted
H\"older continuity. More precisely, we have proven $\p\in C^{0,\mu}_\loc$
and $R\in C^{0,p\mu/2}_\loc$ away from $\Sigma$, and we will improve the
statement about $\p$ in a moment.

The dimension estimate for $\Sigma$ then is a classical result, e.g.\
using \cite[Proposition 9.21]{GM}.
What remains to be proven is the H\"older continuity
of $D\p$. Note that the Euler-Lagrange equation \re{elphi} for $\p$
is a linear elliptic equation with coefficients and right-hand side
depending on $R$. Once we know H\"older continuity of $R$, which we do,
away from $\Sigma$, we are in the realm of classical Schauder estimates
which give us H\"older continuity of $\p$ and even $D\p$
wherever $R$ is H\"older
continuous. A version of that fact that fits our need precisely is
\cite[Theorem 5.19]{GM} which reads as follows.
{\em Let $u\in W^{1,2}_\loc(\Omega,\R^m)$ be a solution to
  \[ \d_\a(A^{\a\b}_{ij}\d_\b u^j)=-\d_\a F^\a_i, \]
  with $A^{\a\b}_{ij}\in C^{0,\sigma}_\loc(\Omega)$ satisfying the
  Legendre-Hadamard condition, for some $\sigma\in(0,1)$. If $F^\a_i
  \in C^{0,\sigma}_\loc(\Omega)$, then we have $Du\in C^{0,\sigma}_\loc(\Omega)$.}

This applies to $\p$ and all $\sigma<\frac2p$, and our theorem is proven.\qed

{\bf Remark. }The case $p>3$ is much simpler, since then $W^{1,p}$ already
embeds into some H\"older space. We have full H\"older regularity of $R$
in $C^{0,p/3-1}_\loc(\Omega)$ then,
hence $\p\in C^{1,p/3-1}_\loc(\Omega)$, and the H\"older exponents can be improved
by arguments from this section. We leave the
details to the interested reader.

\section{Dimension reduction for the singular set}\la{sect:five}

We try to follow part (b) of Luckhaus' theorem, where additional assumptions
allow to prove that the singular set has smaller Hausdorff dimension than
what can be estimated by the arguments of the previous section.
Two things prevent us from applying Luckhaus' theorem directly.
Again, our integrand that has two summands of different homogeneity
is not allowed in Luckhaus' assumptions, and his arguments have to be
modified accordingly. Moreover, for the dimension reduction Luckhaus
has to assume that the unknown functions take their values in a compact
Riemannian manifold, which is not the case for our $\p$. However, the
case of $\p$ being allowed to take values in all of $\R^3$ is sufficiently
easy to be included in Luckhaus' reasoning with only minor modifications.


In order to control in $W^{1,2}\cap W^{1,p}$ some blowup sequence
$(\p_i,R_i)$ we are going to use.
we need a {\em monotonicity formula. }Such formulae have played a central
role in the regularity theory for many functionals. The first monotonicity
formula for harmonic maps has surfaced in the physics literature
\cite{GRSB}. Our monotonicity formula controls
$r^{p-3}\int_{B_r(x_0)}|D\p|^2\,dx$ and
$r^{p-3}\int_{B_r(x_0)}|DR|^p\,dx$ for minimizers $(\p,R)$ of $J$. Note
that $r^{-1}\int_{B_r(x_0)}|D\p|^2\,dx$ would be more natural, due to scaling
invariance. But since we cannot deal with both integrands separately,
we are forced to use the common factor $r^{p-3}$ for both.

\begin{lemma}[monotonicity formula]\la{mon}
  Assume $p\in[2,3)$ and that
  $(\p,R)\in W^{1,2}(\Omega,\R^3)\times W^{1,p}(\Omega,SO(3))$
  minimizes $J$. Then, for every $B_s(x_0)\subset\Omega$ and every $r\in
  (0,s)$, we have
  \eas
  \lll{(s^{p-3}+1)\int_{B_s(x_0)}(|P(R^tD\p)|^2+|DR|^p)\,dx
  -r^{p-3}\int_{B_r(x_0)}(|P(R^tD\p)|^2+|DR|^p)\,dx}\\
  &\ge&\int_{B_s\setminus B_r}|x|^{p-3}\,(|DR|^{p-2}\,|\drad R|^2+Q(\p,R))\,dx
  -cr^{p-1/2}(1+\|D\p\|_{L^2(\Omega)}),
  \eeas
  where here
  \[
  Q(\p,R):=|P(R^tD\p)|^2-|P(R^t(D\p-\tf{x-x_0}{|x-x_0|}\otimes\drad\p))|^2\ge0.
  \]
  If $p=3$, the formula holds without the first ``${}+1$''.
\end{lemma}

{\bf Proof.}
Let $t\in(0,s)$. We abbreviate $B_t:=B_t(x_0)$ and assume $x_0=0$ to
shorten notation.
We compare $(\p,R)$ with $(\p_t,R_t):B_s\to\R^3\times SO(3)$
defined by
\[
  \p_t(x):=\left\{\begin{array}{ll}
    \p(\frac{t}{|x|}\,x) & \mbox{ if }0<|x|<t,\\
    \p(x) & \mbox{ if }t\le|x|<s,
  \end{array}\right.\qquad
  R_t(x):=\left\{\begin{array}{ll}
    R(\frac{t}{|x|}\,x) & \mbox{ if }0<|x|<t,\\
    R(x) & \mbox{ if }t\le|x|<s.
  \end{array}\right.
\]
By $\drad$ we denote the radial derivative in the direction of
$\frac{x}{|x|}$.
Using the fact that $|DR-\frac{x}{|x|}\otimes\drad R|^2
=|DR|^2-|\drad R|^2$, we calculate
\eas
  \lll{\frac{3-p}t\int_{B_t}|DR_t|^p
  =\frac{m-p}t\int_0^t\int_{\d B_\tau}\Big(\Big|DR\Big(\frac{t}{\tau}\,x\Big)\Big|^2
    -\Big|\drad R\Big(\frac{t}{\tau}\,x\Big)\Big|^2\Big)^{p/2}\dH^2(x)\,d\tau}\\
  &\le&\frac{3-p}t\int_0^t\int_{\d B_\tau}\Big|DR\Big(\frac{t}{\tau}\,x\Big)\Big|^{p-2}
  \Big(\Big|DR\Big(\frac{t}{\tau}\,x\Big)\Big|^2
    -\Big|\drad R\Big(\frac{t}{\tau}\,x\Big)\Big|^2\Big)\dH^2(x)\,d\tau\\
  &=&\frac{3-p}t\int_0^t\int_{\d B_t}|DR|^{p-2}(|DR|^2-|\drad R|^2)\dH^2
    \frac{\tau^{2-p}}{t^{2-p}}\,d\tau\\
  &=&\int_{\d B_t}|DR|^{p-2}(|DR|^2-|\drad R|^2)\dH^2\\
  &=&\frac{d}{dt}\int_{B_t}|DR|^p\,dx-\int_{\d B_t}|DR|^{p-2}|\drad R|^2
    \dH^2.
\eeas  
The same way, but using $|P(R^t(D\p-\frac{x}{|x|}\otimes\drad\p))|^2
=|P(R^tD\p)|^2-Q(\p,R)$ this time, we also have
\ea{mm1}
  \frac{3-p}t\int_{B_t}|P(R_t^tD\p_t)|^2\,dx&\le&
  \frac1t\int_{B_t}|P(R_t^tD\p_t)|^2\,dx\nn\\
  &\le&\frac{d}{dt}\int_{B_t}|P(R^tD\p)|^2\,dx-\int_{\d B_t}Q(\p,R)\dH^2.
\eea
We use the fact that $(\p,R)$ minimizes $J$ on $B_t$ and coincides with
$(\p_t,R_t)$ on $\d B_t$. This implies
\ea{mopt}
\lll{\int_{B_t}|P(R^tD\p)|^2\,dx+\int_{B_t}|DR|^p\,dx}\nn\\
  &\le&\int_{B_t}|P(R^tD\p-I)|^2\,dx+\int_{B_t}|P(I)|^2\,dx\nn\\
  &&{}+c\Big(\int_{B_t}|P(R^tD\p-I)|^2\,dx\int_{B_t}|P(I)|^2\,dx\Big)^{1/2}
  +\int_{B_t}|DR|^p\,dx\nn\\
  &\le&\int_{B_t}|P(R^tD\p-I)|^2\,dx+\int_{B_t}|DR|^p\,dx
  +ct^3+ct^{3/2}\|D\p\|_{L^2(\Omega)}\nn\\
  &\le&\int_{B_t}|P(R_t^tD\p_t-I)|^2\,dx+\int_{B_t}|DR_t|^p\,dx
  +ct^3+ct^{3/2}\|D\p\|_{L^2(\Omega)}\nn\\  
  &\le&(1+\eps t)\int_{B_t}|P(R_t^tD\p_t)|^2\,dx+c\eps^{-1}t^2+
  \int_{B_t}|DR_t|^p\,dx\nn\\
  &&{}+ct^{3/2}(1+\|D\p\|_{L^2(\Omega)})\nn\\
  &\le&\int_{B_t}|P(R_t^tD\p_t)|^2\,dx+\int_{B_t}|DR_t|^p\,dx
  +t\,\frac{d}{dt}\int_{B_t}|P(R^tD\p)|^2\,dx\nn\\
  &&{}+ct^{3/2}(1+\|D\p\|_{L^2(\Omega)}),
\eea
where we have chosen $\eps>0$ small enough and used \re{mm1} in the last step.
Now we combine the inequalities and get
\eas
\lll{\frac{d}{dt}\Big(t^{p-3}\int_{B_t}(|P(R^tD\p)|^2+|DR|^p)\,dx\Big)}\\
&=&t^{p-3}\Big(\frac{d}{dt}\int_{B_t}(|P(R^tD\p)|^2+|DR|^p)\,dx
-\frac{3-p}t\int_{B_t}(|P(R^tD\p)|^2+|DR|^p)\,dx\Big)\\
&\ge&t^{p-3}\Big(\frac{d}{dt}\int_{B_t}(|P(R^tD\p)|^2+|DR|^p)\,dx
-\frac{3-p}t\int_{B_t}(|P(R_t^tD\p_t)|^2+|DR_t|^p)\,dx\Big)\\
&&\qquad{}-t^{p-2}\,\frac{d}{dt}\int_{B_t}|P(R^tD\p)|^2\,dx
-ct^{p-3/2}(1+\|D\p\|_{L^2(\Omega)})\\
&\ge&t^{p-3}\int_{\d B_t}(|DR|^{p-2}|\drad R|^2+Q(\p,R))\dH^2\\
&&\qquad{}-\frac{d}{dt}\int_{B_t}|P(R^tD\p)|^2\,dx
-ct^{p-3/2}(1+\|D\p\|_{L^2(\Omega)})
\eeas
Integrating from $r$ to $s$, we infer
\eas
  \lll{s^{p-3}\int_{B_s}(|P(R^tD\p)|^2+|DR|^p)\,dx
  -r^{p-3}\int_{B_r}(|P(R^tD\p)|^2+|DR|^p)\,dx}\nn\\
  &\ge&\int_{B_s\setminus B_r}|x|^{p-3}\,(|DR|^{p-2}\,|\drad R|^2+Q(\p,R))\,dx
  -\int_{B_s\setminus B_r}|P(R^tD\p)|^2\,dx\\
  &&\qquad{}-cs^{p-1/2}(1+\|D\p\|_{L^2(\Omega)})
\eeas
for $0<r<s<\dist(x_0,\d\Omega)$. This proves the lemma.\qed


The following proposition summarizes what we can infer from the monotonicity
formula via a blow-up argument. We define the singular set $\Sing(\p,R)$ as
the set of points in $\Omega$ at which $(\p,R)$ fails to be locally
in $C^{1,\mu}\times C^{0,\mu}$ for any $\mu\in(0,1)$.

\begin{proposition}[partial regularity and $p$-minimizing tangent
    maps]\la{hoel2}
  Assume $p\in[2,3]$ and that $(\p,R)\in W^{1,2}(\Omega,\R^3)\times
  W^{1,p}(\Omega,SO(3))$ is a minimizer of $J$ on $\Omega$.
  Then the following statements hold.
  
  (i) The singular set $\Sing(\p,R)$ is discrete in $\Omega$. If $p=3$,
  it is even empty.

  (ii) If $p\in(2,3)$ and $\Sing(\p,R)$ is not empty, then there is
  at least one ``$p$-minimizing tangent map'' to $SO(3)$, by which we mean
  a nonconstant continuous map $R_\infty:B^3\setminus\{0\}\to SO(3)$ which is
  radially constant and minimizes $E^p(\tR):=\int_{B^3}|D\tR|^p\,dx$
  among all $W^{1,p}$-maps to $SO(3)$ with the same boundary values.

  (iii) The same statement holds if $p=2$ and $\mue=\muc=\muz$.
\end{proposition}

The notion of $p$-minimizing tangent maps has been central in the regularity
for $p$-harmonic maps. Like here, they are the obstacles to full regularity
of $p$-harmonic maps and hence play an important role in \cite{Lu} and
in many other results on $p$-harmonic maps. And, of course, $p$-minimizing
tangent maps are weakly $p$-harmonic maps themselves.

{\bf Proof.}
By Proposition \ref{hoel1}, the singular set is a subset of the set
$\Sigma$ where ``enough $p$-energy of $R$ concentrates''. In what follows,
we assume that $x_0\in\Sigma$. We denote rescaled versions of $\p$
and $R$ around $x_0$ by
\[
  \p_i(x):=\rho_i^{p/2-1}(\p(x_0+\rho_i x)-\opi),\qquad
  R_i(x):=R(x_0+\rho_i x),
\]
where $(\rho_i)_{i\in\N}$ is any strictly decreasing sequence with
$\rho_1$ sufficiently small and $\rho_i\searrow0$.
Then $(\p_i,R_i)$ minimizes
\[
J_i(\tp,\tR):=\int_{B^3}(|P(\tR^tD\tp-\rho_i^{p/2}I)|^2+|D\tR|^p)\,dx,
\]
and this explains the scaling chosen, because the monotonicity formula
from Lemma \ref{mon} implies that $J_i(\p_i,R_i)$ stays bounded as
$i\to\infty$.

Assuming $\rho_i\searrow0$, the formal limit functional is
\[
J_\infty(\tp,\tR):=\int_{B^3}(|P(\tR^tD\tp)|^2+|D\tR|^p)\,dx.
\]
Note that this time all $R_i$ map to $SO(3)$ rather than to rescaled
and shifted copies of it. The bound on $J_i(\p_i,R_i)$ and the fact that
$\p_i$ has mean value $0$ imply that $(\p_i,R_i)$ is bounded in
$W^{1,2}(B^3,\R^3)\times W^{1,p}(B^3,SO(3))$. After passing to a subsequence,
we can therefore assume that $(\p_i,R_i)\rightharpoonup(\p_\infty,R_\infty)$
weakly in $W^{1,2}\times W^{1,p}$ for some $(\p_\infty,R_\infty)\in
W^{1,2}(B^3,\R^3)\times W^{1,p}(B^3,SO(3))$.

By the same reasoning as in Lemma \ref{comp} (which we do not work out
since it is sufficiently close to \cite{Lu}), we have
convergence $(\p_i,R_i)\to(\p_\infty,R_\infty)$ in $W^{1,2}\times W^{1,p}$-norm,
and $(\p_\infty,R_\infty)$ minimizes $J_\infty$ with respect to its boundary
values.

Inserting $s=\rho_i$, $r=\rho_j$ into Lemma \ref{mon} and rescaling, we have
\eas
  \lll{(1+\rho_i^{3-p})\int_{B_1}(|P(R^tD\p_i)|^2+|DR_i|^p)\,dx
  -\int_{B_1}(|P(R^tD\p_j)|^2+|DR_j|^p)\,dx}\\
  &\ge&\int_{B_1\setminus B_{\rho_j/\rho_i}}|x|^{p-3}\,(|DR_i|^{p-2}\,|\drad R_i|^2
  +Q(\p_i,R_i))\,dx
  -c\rho_j^{p-1/2}
  \qquad\qquad\strut
\eeas
for $2\le p<3$, and the same with $(1+\rho_i^{3-p})$ replaced by $1$ for
$p=3$.
Letting $j\to\infty$ and then $i\to\infty$, the norm convergence gives
\[
\int_{B^3}|x|^{p-3}\,(|DR_\infty|^{p-2}\,|\drad R_\infty|^2
+Q(\p_\infty,R_\infty))\,dx=0.
\]
Hence $\drad R_\infty\equiv 0$, and $Q(\p_\infty,R_\infty)\equiv0$, which also
implies $\drad\p_\infty\equiv0$. This means that both $\p_\infty$ and
$R_\infty$ are radially constant. And by the definition of $\Sigma$ and
the norm convergence, $(\p_\infty,R_\infty)$ is not constant.

Our first aim is to prove that the singular set of $(\p,R)$ is discrete.
The strategy for that is taken from \cite{Lu}, but goes back to
``Federer's dimension reduction argument'' from \cite{Fe}.
This is clearly a local question that we can answer by considering
only a neighborhood of some point in $\Omega$. Therefore we restrict
to $\Omega=B^3$ and define the ``$\eps_0$-singular set'' of any
$(\p,R)$ by
\eas
&&S(\p,R;\eps_0)\;:=\;\Big\{x_0\in B^3:|x_0|<\frac12,\\
&&\quad r^{-1}\|D\p\|_{L^2(B_\rho(x_0))}^2+
r^{p-3}\|DR\|_{L^p(B_\rho(x_0))}^p\ge\eps_0\mbox{ for }0<r<1-|x_0|\Big\}\,.
\eeas
By Proposition \ref{hoel1}, it is sufficient to show that $S(\p,R;\eps_0)$
is discrete for $\eps_0>0$ chosen sufficiently small.

Suppose the contrary, then there is an $a\in S(\p,R;\eps_0)$ that is the limit
of a sequence $\{a_i\}_{i\in\N}$ in $S(\p,R;\eps_0)\setminus\{a\}$.
Then define rescaled mappings $(\p_i,R_i)$ as above (with $x_0$ replaced by
$a$), for some $\rho_i\searrow0$ that the sequence converges
to a radially constant minimizer $(\p_\infty,R_\infty)$.
Since the integrals transform naturally
under the scaling involved in the definition of $\p_i$ and $R_i$, we find that
$|x_0|\le\frac12$ and $a+\rho_i^{-1}x_0\in S(\p,R;\eps_0)$ imply
$x_0\in S(\p_i,R_i;\eps_0)$.

And we can also assume that we have chosen
the $\rho_i$ in such a way that the sequence $\rho_i(a_i-a)$ has an
accumulation point $x_a$ with $|x_a|=\frac12$. Then we can verify that
$x_a\in S(\p_\infty,R_\infty;\eps_0)$.
But $x_a\ne0$, and $(\p_\infty,R_\infty)$ is radially constant, and
nonconstant because of the definition of $S(\p,R,\eps)$ and the
norm convergence. Therefore,
we have $\R x_a\subseteq S(\p_\infty,R_\infty;\eps_0)$, but the $\eps_0$-singular
set of a $W^{1,p}$-mapping for $p\ge2$ in three dimensions always has
vanishing one-dimensional Hausdorff measure. This is a contradiction,
and we have proved Assertion (i) that the $\epsilon_0$-singular set of
$(\p,R)$, and hence also the singular set, is discrete. And as
the $\eps_0$-singular set of a $W^{1,3}$-mapping always vanishes in three
dimension, the singular set is empty in case $p=3$.


Now let us assume $p\in(2,3)$ and prove (ii). To this end, we remark
that $(\p_\infty,R_\infty)$ minimizes $J_\infty$, a functional very similar
to $J$ except for being less inhomogeneous for the lack of the term
involving $I$. In particular, a monotonicity formula for minimizers
of $J_\infty$ is proven exactly as for those of $J$, the proof being
actually slightly shorter since there is one term less involved.
That monotonicity formula allows us to blow up $(\p_\infty,R_\infty)$
again around $0$, finding that, for every $i\in\N$, the pair
\[
\p_{\infty,i}(x):=\rho_i^{p/2-1}\p_\infty(\rho_i x),\qquad
R_{\infty,i}(x):=R_\infty(\rho_i x)
\]
minimizes $J_\infty$ with respect to its boundary values. The compactness
argument employed above gives that the weak limit again minimizes $J_\infty$,
but by the radial homogeneity of $\p_\infty$ and $R_\infty$, we have
\[
(\p_{\infty,i},R_{\infty,i})\to(0,R_\infty),
\]
This means that $(0,R_\infty)$ is a $J_\infty$-minimizer, which of course implies
that $R_\infty$ minimizes $\int_{B^3}|D\tR|^p\,dx$ and hence is a weakly
$p$-harmonic map. By its properties already known, including $R_\infty$ being
nonconstant by repetition of the argument above, it is a $p$-minimizing
tangent map, and we have proved (ii).

The argument we just performed breaks down for $p=2$, since this time
we have actually taken advantage of the inhomogeneity of the functional.
We do not know whether the same statement can be proven for $p=2$, except
for the special case when $\mue=\muc=\muz$. In that case,
the functional $J_\infty$ equals
\[
J_\infty(\tp,\tR)=\int_{B^3}(\mue|D\tp|^2+|D\tR|^p)\,dx,
\]
in which $\tp$ and $\tR$ are ``decoupled''. Which means that
$(\p_\infty,R_\infty)$ minimizes $J_\infty$ iff $\p_\infty$ minimizes
$\|D\tp\|_{L^2(B^3)}^2$ and $R_\infty$ minimizes $\|D\tR\|_{L^p(B^3)}^p$.
This allows the same conclusion as in (ii) and hence proves (iii),
which completes the proof of Proposition \ref{hoel2}.\qed

\section{Nonexistence of $p$-harmonic minimizing tangent maps}
\la{sect:six}

We have seen in Proposition \ref{hoel2} that, in order to exclude point
singularities, we need to know the non-existence of $p$-minimizing
tangent maps (now abbreviated as $p$-mtm)
$B^3\setminus\{0\}\to SO(3)$. The first thing we claim is that it
suffices to exclude $p$-mtm to $S^3$ instead.

\begin{lemma}[lift of $p$-mtm to $S^3$]\la{lift}
  Assume $p>1$, and that there is a $p$-mtm $R:B^3\setminus\{0\}
  \to SO(3)$. Then there also is a $p$-mtm $u:B^3\setminus\{0\}\to S^3$.
\end{lemma}

{\bf Proof.} This is due to the fact that $S^3$ is the universal cover
of $SO(3)$. A locally isometric (up to a constant factor) $2$-to-$1$
covering map $\pi:S^3\to SO(3)$ is given by
\[
  \pi(w,x,y,z):=\left(\begin{array}{ccc}
    1-2y^2-2z^2 & 2xy-2zw & 2xz+2yw \\
    2xy+2zw & 1-2x^2-2z^2 & 2yz-2xw \\
    2xz-2yw & 2yz+2xw & 1-2x^2-2y^2
  \end{array}\right)\,.
\]
Since $B^3\setminus\{0\}$ is simply connected, there is a lift
$u:B^3\setminus\{0\}\to S^3$ of our given map $R$, satisfying
$\pi\circ u=R$. Certainly, $u$ cannot be constant if $R$ is not.
And $u$ minimizes $E^p$ with respect to its own boundary values on $S^2$.
Assume it does not, then we had $v:B^3\to S^3$ with $v=u$ on $S^2$ which
has $E^p(v)<E^p(u)$. But then, since $\pi$ is a constant times an isometry,
$E^p(\pi\circ v)<E^p(\pi\circ u)=E^p(R)$ with $\pi\circ v=R$ on $S^2$,
which means $R$ would not be minimizing. This is a contradiction, hence
$u$ minimizes $E^p$. It is also radially constant. Summarizing, we have
found that $u$ is a $p$-mtm.\qed

{\bf Remark.} The proof gives some interpretation of the singular
example for
the Cosserat model given in Section \ref{sect:ex}. The mapping $R$ there is
just the equator map $(\frac{x}{|x|},0):B^3\to S^3$ projected to
$SO(3)$ using $\pi$.\qed

The Lemma means that we know H\"older continuity of minimizers in our
Cosserat model once we can exclude the existence of $p$-mtm
$B^3\setminus\{0\}\to S^3$. But this has been done for $p=2$.

\begin{proposition}[a complete regularity case for $p=2$]\la{reg1}
  Assume $p=2$, $\mue=\muc=\muz$, and that $(\p,R)\in
  W^{1,2}(\Omega,\R^3\times SO(3))$ is a minimizer of $J$ on $\Omega$.
  Then $(\p,R)\in C^{1,\mu}_\loc(\Omega,\R^3)\times C^{0,\mu}_\loc(\Omega,SO(3))$
  for every $\mu\in(0,1)$.
\end{proposition}

{\bf Proof.} This is a combination of Proposition~\ref{hoel1},
Proposition~\ref{hoel2} and the nonexistence of $2$-mtm
$B^3\setminus\{0\}\to S^3$. The latter has been proven by Schoen
and Uhlenbeck in \cite[Theorem~2.7]{SU2}.\qed

In order to apply part (ii) of Proposition \ref{hoel2}, we would like
to know about nonexistence of $p$-mtm $B^3\setminus\{0\}\to S^3$.
Unfortunately, the proof of \cite{SU2} does not seem to give a hint
on how to handle that, because it uses the fact that harmonic $2$-spheres
are very closely related to minimal immersions, an argument that does
not carry over to $p$-harmonic maps for any $p\ne2$. Okayasu \cite{Ok}
has published a modification of the \cite{SU2}-argument for the target
$S^3$ (and others), based on a so-called improved Kato inequality. He
proves that, for harmonic maps $S^{k-1}\to N$, the trivial pointwise inequality
$|D|Du||\le|\nabla Du|$ can be improved to
$|D|Du||^2\le\frac{k-2}{k-1}\,|\nabla Du|^2$,
and the constant is optimal.

If we had an improved Kato-type inequality for $p$-harmonic maps, we could
try to improve regularity theorems from Xin and Yang \cite{XY} or
Nakauchi \cite{Na} the way Okayasu improved \cite{SU2}. Unfortunately,
no optimal Kato inequality for $p$-harmonic maps seems to be available.
There is, however, an optimal one for $p$-harmonic {\em functions\/}
(i.e.\ maps to $\R$) that can give is some orientation. It has
been proven recently by Chang, Chen, and Wei
\cite[Lemma 5.4]{CCW} and reads $|\nabla Du|^2\ge(1+\tilde{\kappa})|D|Du||^2$,
where here $\tilde{\kappa}:=\min\{\frac{(p-1)^2}{m-1}\,,1\}$, and $m$ is the
domain dimension. We do not see
how the proof could carry over to the $S^3$-valued case, but we do get some
improvement of Kato's inequality, which is probably not optimal, but for
$p\searrow2$ reproduces Okayasu's result.

\begin{lemma}[improved Kato inequality for $p$-harmonic maps]\la{kato}
  Let $p>1$, and  let $M$ and $N$ be smooth complete Riemannian manifolds,
  $m:=\dim M$. Fix some parameter $\eps\in(0,1]$, let 
  \[
    \kappa:=\frac{m-1+(\frac1{\eps}-1)(p-2)^2}{m-\eps}
    \qquad\mbox{(which is $>\frac12$)},
  \]
  and assume that $u\in C^{1,\a}(M,N)$ is a $p$-harmonic mapping
  (which is automatically $C^2$ away from the points with
  $Du(x)=0$). Then at any $x\in M$ with $du(x)\ne0$, we have
  \[
  |D|Du||^2\le\kappa\,|\nabla Du|^2,
  \]
  where $\nabla$ denotes the Levi-Civita connection of $N$.
\end{lemma}

{\bf Proof. } 
By Nash's embedding theorem, $N$ (or some compact portion of it
around $u(x)$) is embedded isometrically into some $\R^n$. We write
$\nabla_i$ for partial derivatives with respect to the Levi-Civita connection
of $N$.
The $p$-harmonic map equation can be written as
\[
\sum_i\nabla_i(|Du|^{p-2}\d_iu)=0,
\]
which is equivalent to
\[
\sum_i\nabla_i\d_i u=\frac{p-2}{|Du|^2}\,\sum_{i,j}\<\nabla_i\d_ju,\d_ju\>\d_iu.
\]
This implies, abbreviating $\sum_i\nabla_i\d_i=:\tau$, known as the
{\em tension field\/},
\eq{k1}
|\tau(u)|^2\le(p-2)^2|\nabla Du|^2.
\eeq
We now fix $x\in M$ and some index $\a\in\{1,\ldots,n\}$ and find an ONB
$\{b_1,\ldots, b_n\}$ of $T_xM$ (depending on $\a$) such that
$\d_1u^\a=|Du^\a|b_1$, which implies $\d_1u^\a(x)=|Du^\a(x)|$ and $\d_ju^\a(x)=0$
for $j\ne1$. This idea is from \cite{CCW}. By using appropriate coordinates,
we can also assume $\nabla_i\d_ju^\a(x)=\nabla_j\d_iu^\a(x)$ for all $i,j$.
\eas
\lll{\sum_{i,j=1}^m(\nabla_i\d_ju^\a)^2}\\
&\ge&(\nabla_1\d_1u)^2+2\sum_{h=2}^m(\nabla_h\d_1u^\a)^2
  +\sum_{h=2}^m(\nabla_h\d_hu^\a)^2\\
&\ge&(\nabla_1\d_1u)^2+2\sum_{h=2}^m(\nabla_h\d_1u^\a)^2
  +\frac1{m-1}\,\Big(\sum_{h=2}^m\nabla_h\d_hu^\a\Big)^2\\
&=&(\nabla_1\d_1u)^2+2\sum_{h=2}^m(\nabla_h\d_1u^\a)^2
  +\frac1{m-1}\,(\tau(u)^\a-\nabla_1\d_1u^\a)^2\\
&\ge&\frac{m-\eps}{m-1}\,(\nabla_1\d_1u^\a)^2+2\sum_{h=2}^m(\nabla_h\d_1u^\a)^2
  -\frac{\eps^{-1}-1}{m-1}\,(\tau(u)^\a)^2.
\eeas
In the last line, we have applied Young's inequality. Using $2\ge
\frac{m-\eps}{m-1}$ and $(\nabla_h\d_1u^\a)^2=(\d_h|Du^\a|)^2$,
the last estimate becomes
\eq{k2}
\sum_{i,j=1}^m(\nabla_i\d_ju^\a)^2\ge\frac{m-\eps}{m-1}|D|Du^\a||^2
-\frac{\eps^{-1}-1}{m-1}\,(\tau(u)^\a)^2.
\eeq
We have
\eas
|D|Du||^2&=&\Big|D\sqrt{\sum_\a|Du^\a|^2}\Big|^2
\;=\;\Big(\sum_\a\frac{|D|Du^\a||\,|Du^\a|}{|Du|}
\Big)^2\\&\le&\sum_\a|D|Du^\a||^2\sum_\a\frac{|Du^\a|^2}{|Du|^2}
\;=\;\sum_\a|D|Du^\a||^2
\eeas
and can therefore sum over $\a$ in \re{k2}. Using also \re{k1}, we get
\eas
|\nabla Du|^2&\ge&\frac{m-\eps}{m-1}\,|D|Du||^2
-\frac{\eps^{-1}-1}{m-1}\,|\tau(u)|^2\\
&\ge&\frac{m-\eps}{m-1}\,|D|Du||^2
-\frac{\eps^{-1}-1}{m-1}\,(p-2)^2|\nabla Du|^2.
\eeas
Absorbing the last term into the left-hand side proves the lemma.\qed

{\bf Remark. }
If $p-2$ is moderately large, there is no choice for $\eps>0$ that makes
$\kappa<1$. Therefore, Lemma \ref{kato} must be seen as a tool
only for $p$ close to $2$. If $p>2$ is close to $2$, the optimal choice
of $\eps$, the one that makes $\kappa$ smallest, is
\[
\eps:=\frac{1}{m-1-(p-2)^2}\Big((p-2)^2
+\sqrt{(p-2)^4+m(m-1-(p-2)^2)(p-2)^2}\Big)\,.
\]
For $p\searrow2$, we have $\eps\sim\sqrt{\frac{m}{m-1}}\,(p-2)$.\qed

Now we use a slight modification of the arguments in \cite{XY} or \cite{Na}
--- the latter is slightly easier to cite for our purposes. We assume
that $u:B^k\setminus\{0\}\to S^n$ is a minimizing $p$-harmonic tangent map.
Then we use \cite[Lemma 1]{Na}, which is derived from the stability
inequality,
\eq{na1}
\int_{S^{k-1}}|Du|^{p-2}|D|Du||^2\,dx
\ge\frac{n-p}{n+p-2}\int_{S^{k-1}}|Du|^{p+2}\,dx
-\frac{(k-p)^2}4\int_{S^{k-1}}|Du|^p\,dx.
\eeq
The other ingredient is \cite[Lemma 2]{Na}, which we cite in a modified version.
We use Nakauchi's derivation of the lemma, but do not use 
$|\nabla Du|\ge|D|Du||$, and find
\eas
\lll{\int_{S^{k-1}}|Du|^{p-2}(|\nabla Du|^2+(p-2)|D|Du||^2)\,dx}\\
&\le&\frac{k-2}{k-1}\int_{S^{k-1}}|Du|^{p+2}\,dx
-(k-2)\int_{S^{k-1}}|Du|^p\,dx.
\eeas
Using the improved Kato inequality from Lemma \ref{kato}, we find
\ea{na2}
\lll{\Big(\frac{k-1-\eps}{k-2+(\eps^{-1}-1)(p-2)^2}+p-2\Big)
\int_{S^{k-1}}|Du|^{p-2}|D|Du||^2\,dx}\nn\\
&\le&\frac{k-2}{k-1}\int_{S^{k-1}}|Du|^{p+2}\,dx
-(k-2)\int_{S^{k-1}}|Du|^p\,dx.
\eea
Comparing \re{na2} with \re{na1}, we could find slight improvements of the
results of \cite{XY} and \cite{Na}, but we do not bother to write them down
in general form. Instead, we concentrate on the case we need to exclude
point singularities in the Cosserat model, and for that we need only consider
the special case $k=n=3$. Combining \re{na1} and \re{na2}, we then have
\eas
\lll{\Big(\frac{3-p}{1+p}\Big(\frac{2-\eps}{1+(\eps^{-1}-1)(p-2)^2}+p-2\Big)
-\frac12\Big)\int_{S^2}|Du|^{p+2}\,dx}\\
&\le&\Big(\frac{(3-p)^2}4\Big(\frac{2-\eps}{1+(\eps^{-1}-1)(p-2)^2}+p-2\Big)
-1\Big)\int_{S^2}|Du|^p\,dx.
\eeas
If the coefficient on the left-hand side is $>0$ while the one on the
right-hand side is $\le0$, then $u$ must be constant, hence there exist no
$p$-mtm. (This is the key idea from \cite{SU2}.)

Using the optimal $\eps$ from the previous remark,
this condition is easily verified numerically for $2\le p\le\frac{32}{15}$,
where the upper bound is not quite optimal and set to a fraction
only for simplicity --- we do not expect our method of proof to be
optimal, anyway.

This means that for $2\le p\le\frac{32}{15}$, there are no nonconstant
$p$-minimizing tangent maps $B^3\setminus\{0\}\to S^3$.
Combining that with Proposition \ref{hoel1}, Proposition \ref{hoel2} (ii),
and Lemma \ref{lift}, we have proven the following result.

\begin{proposition}[H\"older regularity for $p>2$ close to $2$]\la{reg2}
  Assume $p\in(2,\frac{32}{15}]$ and that $(\p,R)\in
    W^{1,2}(\Omega,\R^3)\times W^{1,p}(\Omega,SO(3))$ is a minimizer of $J$ on
    $\Omega$.
    Then $(\p,R)\in C^{1,\mu}_\loc(\Omega,\R^3)\times C^{0,\mu}_\loc(\Omega,SO(3))$
  for every $\mu\in(0,\frac2p)$.\hfill$\Box$
\end{proposition}

{\bf Remark. }From the Propositions \ref{reg1} and \ref{reg2}, we read off
that the singular example from Section \ref{sect:ex} cannot minimize
$J$ if $p\in[2,\frac{32}{15}]$.\qed

\section{External forces and moments}\la{sect:forces}

Now we return to the case where the functions $f$ and $M$ do not vanish.
That is, we again consider the functional
\[
  J(\p,R):=\int_\Omega\Big(|P(R^tD\p-I)|^2+|DR|^p
  +\p\cdot f+R\cdot M\Big)\,dx,
\]
where we have omitted a $-\int x\cdot f$ which only gives an additive
constant and is therefore irrelevant for minimizing.

Then we have to care for some more lower order terms which do not
really affect the reasoning of the previous sections. First of all,
the monotonicity allows for force and moment terms in the natural
Lebesgue spaces the functional allows.

\begin{lemma}[more general monotonicity formula]
  Additionally to the assumptions of Lemma \ref{mon}, let functions
  $f\in L^2(\Omega,\R^3)$ and $M\in L^q(\Omega,\R^{3\times3})$ be given,
  where here $\frac1p+\frac1q=1$.
  Let $(\p,R)$ be a minimizer of the functional $J$ now involving
  the corresponding force and moment potentials. Then the monotonicity
  formula from Lemma \ref{mon} still holds after it has been modified
  by an additional term
  \[
  -cs^{p-1}(\|D\p\|_{L^2(\Omega)}\|f\|_{L^2(\Omega)}
  +\|DR\|_{L^p(\Omega)}\|M\|_{L^q(\Omega)})
  \]
  on the right-hand side.
\end{lemma}

{\bf Proof.} The change to be made in the proof of Lemma \ref{mon} is
in the estimate \re{mopt} where $(\p,R)$ is compared with $(\p_t,R_t)$.
Here we have to add a term $\int_{B_t}(\p_t-\p)\cdot f+(R_t-R)\cdot M$
to the right-hand side, which in the sequel is estimated according to
\eas
\lll{\int_{B_t}((\p_t-\p)\cdot f+(R_t-R)\cdot M)\,dx}\\
&\le&\|\p_t-\p\|_{L^2(B_t)}\|f\|_{L^2(B_t)}+\|R_t-R\|_{L^p(B_t)}\|M\|_{L^q(B_t)}\\
&\le&ct(\|D\p\|_{L^2(\Omega)}\|f\|_{L^2(\Omega)}
+\|DR\|_{L^p(\Omega)}\|M\|_{L^q(\Omega)}).
\eeas
This term is easily carried through the remaining estimates in
the proof of the monotonicity formula.\qed

For the regularity theory, we restrict to H\"older
continuous data, which should be
good enough for most applications. In comparison to the case without
exterior forces, we have to compromise about the H\"older exponents.

\begin{proposition}[regularity with exterior forces and moments]\la{reg3}
  The statements of the Propositions \ref{hoel1}, \ref{hoel2}, \ref{reg1},
  and \ref{reg2} continue to hold in the case of nonvanishing $f$ and $M$
  if we assume $f\in C^{0,\mu}(\overline{\Omega},\R^3)$ and
  $M\in C^0(\overline{\Omega},\R^{3\times3})$. However, if $f\not\equiv0$,
  we have to restrict the H\"older exponent $\mu$
  to $(0,\frac1{2p})$ instead of $(0,\frac{p}2)$.
\end{proposition}

{\bf Sketch of proof.} We have to check the arguments of the proofs
wherever we have used the minimality of $(\p,R)$ or the Euler-Lagrange
equations. As a matter of fact, this does not affect our reasoning too
much, since the proofs use blowup procedures, and the new potential
terms are scaled away in the blowup processes. Hence the limit
functionals are the same as in the case where $f$ and $M$ vanish.

To see how that works, let us first have a look into the proof of
Lemma~\ref{morr}. The functional that $(\p_i,R_i)$ minimizes now
is a modified version of the $J_i$ given there, namely
\eas
J_i(\tp,\tR)&:=&\frac{\g_i^2}{\g_i^2+\del_i^p}\int_{B^3}|P((\oRi+\del_i\tR)^t
D\tp-\rho_i\g_i^{-1}I)|^2\,dx+\frac{\del_i^p}{\g_i^2+\del_i^p}\int_{B^3}|D\tR|^p
\,dx\\
&&{}+\frac{\rho_i^2}{\g_i^2+\del_i^p}\int_{B^3}(\opi+\g_i\tp(x))\cdot f(x_i+
\rho_i x))\dx\\
&&{}+\frac{\rho_i^2}{\g_i^2+\del_i^p}\int_{B^3}(\oRi+\del_i\tR(x))\cdot M(x_i
+\rho_i x))\,dx.
\eeas
Remember $\rho_i^{2\mu}\le\g_i^2+\del_i^p$, and $\mu<1$, hence the coefficient
$\frac{\rho_i^2}{\g_i^2+\del_i^p}$ of the last two integrals vanishes in
the limit $i\to\infty$. And $(\p_\infty,R_\infty)$ minimizes the original
$J_\infty$ without any additional terms once we can prove that the last
two integrals in $J_i(\p_i,R_i)$ are bounded uniformly in $i$. But they are,
because first of all, we have
\[
\Big|\int_{B^3}((\oRi+\del_i\tR(x))\cdot M(x_i+\rho_i x))\,dx\Big|
=\Big|\rho_i^{-3}\int_{B_{\rho_i}(x_i)}R\cdot M\,dx\Big|
\le c\|M\|_{C^0(\overline{\Omega})}.
\]
For the other integral, things are slightly more involved, and we need
$\mu<\frac14$ in Lemma~\ref{morr}, which corresponds to $\mu<\frac1{2p}$
in Proposition~\ref{hoel1}. If $\mu<\frac14$, we have
$\frac{\rho_i^{1/2}}{\g_i^2+\del_i^p}\to0$, and it is sufficient to bound
$\rho_i^{3/2}$ times the integral. We have
\eas
  \lll{\rho_i^{3/2}\Big|\int_{B^3}(\opi+\g_i\p_i(x))\cdot f(x_i+\rho_i x))
    \,dx\Big|}\\
  &=&\rho_i^{-3/2}\Big|\int_{B_{\rho_i}(x_i)}\p\cdot f\,dx\Big|
  \;\le\;c\|\p\|_{L^2(B_{\rho_i}(x_i))}\|f\|_{C^0(\overline{\Omega})}
  \;\le\;c\|\p\|_{L^2(\Omega)}\|f\|_{C^0(\overline{\Omega})}.
\eeas
Note that the estimate depends on $\|\p\|_{L^2(\Omega)}$ for our fixed minimizer
$(\p,R)$. But it is used in a term that vanishes in the limit, anyway, and
we can check the proof of Lemma~\ref{morr} to find that the constants  in its
statement still do not depend on $\|\p\|_{L^2(\Omega)}$.

In the blow-up performed in the proof of Proposition~\ref{hoel2}, there is a
similar reasoning. Here, two integrals have to be added to the functional
$J_i$ that correspond to $\rho_i^{3-p}\int_{B_r(x_0)}\p\cdot f\,dx$ and
$\rho_i^{3-p}\int_{B_r(x_0)}R\cdot M\,dx$, both of which vanish as $i\to\infty$,
even if $p=3$.

There is one more change necessary in the proof of Proposition~\ref{hoel1},
where we apply Schauder theory in order to get H\"older continuity
of $D\p$. The equation considered for this now reads
$\div(RP^2(R^tD\p-I))=f$, with $f\in C^{0,\mu}$ instead of $0$. But this
additional right-hand side in $C^{0,\mu}$ does not affect the Schauder
estimates, which are known to hold also for $f_i-\d_\a F^\a_i$ instead of
just $-\d_\a F^\a_i$ if the $f_i$ and $F^\a_i$ are in $C^{0,\mu}$.

Apart from that, there are no serious changes in the proofs, which can
therefore be adapted to prove Proposition~\ref{reg3}.\qed

\vfill

\begin{center}
\scriptsize
Fakult\"at f\"ur Mathematik, Universit\"at Duisburg-Essen, D-45117 Essen, 
Germany.\\
{\tt andreas.gastel@uni-due.de}
\end{center}

\end{document}